\documentclass[12pt]{elsarticle}
\usepackage{amsmath,amssymb}
\usepackage{eepic}
\usepackage{graphicx}
\usepackage{enumitem}
\biboptions{sort&compress}
\usepackage[left=2.0cm,right=2.0cm,top=1cm,bottom=2cm,bindingoffset=0cm]{geometry}

\journal{arXiv}

\begin{document}

\begin{frontmatter}

\title{
  Hidden and self-excited attractors \\ in electromechanical systems
  with and without equilibria.
}

\author{M.A. Kiseleva}
\author{N.V. Kuznetsov\corref{cor}}
\ead{nkuznetsov239@gmail.com}
\author{G.A.~Leonov}

\address[spbu]{Faculty of Mathematics and Mechanics,
St.~Petersburg State University, Russia}
\address[fin]{Department of Mathematical Information Technology,
University of Jyv\"{a}skyl\"{a}, Finland}
\address[ras]{
Institute for Problems in Mechanical Engineering of the Russian Academy of Sciences,
Russia}

\begin{abstract}
  This paper studies hidden oscillations appearing in electromechanical systems with and without equilibria. Three different systems with such effects are considered: translational oscillator-rotational actuator, drilling system actuated by a DC-motor and drilling system actuated by induction motor. We demonstrate that all three systems experience hidden oscillations in sense of mathematical definition. But from physical point of view in certain cases it is quite easy to localize there oscillations.
\end{abstract}


\end{frontmatter}

\section{Introduction}
An oscillation in a dynamical system is either self-excited or hidden.
Study of stability and oscillations in electromechanical systems requires construction of mathematical model and its analysis. Stability corresponds to normal operation of the system, oscillations are localized in case the initial conditions from their open neighbourhood lead to long-time behaviour that approaches the oscillation.
Depending on simplicity of finding the basin of attraction in the phase space it is natural to suggest the following classification of attractors \citep{KuznetsovLV-2010-IFAC,LeonovKV-2011-PLA,LeonovKV-2012-PhysD,LeonovK-2013-IJBC,LeonovKM-2015-EPJST}:
{\it An attractor is called a \emph{hidden attractor} if its
 basin of attraction does not intersect with
 small neighborhoods of equilibria,
 otherwise it is called a \emph{self-excited attractor}.
}
Self-excited attractor's basin of attraction
is connected with an unstable equilibrium. Therefore, self-excited attractors
can be localized numerically by the
\emph{standard computational procedure},
in which a trajectory,
which starts from a point of an unstable manifold in a neighbourhood
of an unstable equilibrium, after a transient process
is attracted to the state of oscillation and traces it.
In contrast, hidden attractor's basin of attraction
is not connected with unstable equilibria.
For example, hidden attractors are attractors in
the systems with no equilibria
or with only one stable equilibrium
(a special case of multistable systems and
coexistence of attractors). 
Recent examples of hidden attractors can be found in
\emph{The European Physical Journal Special Topics: Multistability: Uncovering Hidden Attractors}, 2015
(see \cite{ShahzadPAJH-2015-HA,BrezetskyiDK-2015-HA,JafariSN-2015-HA,ZhusubaliyevMCM-2015-HA,SahaSRC-2015-HA,Semenov20151553,FengW-2015-HA,Li20151493,FengPW-2015-HA,Sprott20151409,Pham20151507,VaidyanathanPV-2015-HA,SharmaSPKL-2015-EPJST}).
See also 
\citep{KuznetsovKL-2013-DEDS,LeonovK-2011-DAN,BraginVKL-2011,AndrievskyKLP-2013-IFAC,LeonovKKSZ-2014,KuznetsovL-2014-IFACWC,BianchiKLYY-2015,LeonovKKK-2015-IFAC,Kuznetsov-2016,
KingniJSW-2014,BurkinK-2014-HA,PhamRFF-2014-HA,CafagnaG-2015-woeq,WangSWZ-2015-HA,KuznetsovKMS-2015-HA,Zelinka-2016-HA,Chen-2015-IFAC-HA,BaoHCXY-2015-HA}.

Hidden oscillations appear naturally in systems without equilibria, describing
various mechanical and electromechanical models with rotation.
One of the first examples of such models was described by Arnold Sommerfeld in 1902 \citep{Sommerfeld-1902}.
He studied vibrations caused by a motor driving an unbalanced mass
and discovered the resonance capture (Sommerfeld effect).
The Sommerfeld effect 
represents the failure of a rotating mechanical system
to be spun up by a torque-limited rotor to a desired rotational velocity due to its resonant
interaction with another part of the system \citep{EvanIwanowski-1976,Eckert-2013}.
Relating this phenomenon
to the real world Sommerfeld wrote, ``{\it This experiment corresponds roughly to
the case in which a factory owner has a machine set on a poor foundation running
at 30 horsepower. He achieves an effective level of just 1/3, however, because only 10
horsepower are doing useful work, while 20 horsepower are transferred to the foundational masonry}''
\citep{Eckert-2013}. 
Another well-known chaotic physical system with no equilibrium points
is the Nos\`{e}--Hoover oscillator \cite{Nose-1984,Hoover-1985,Sprott-1994,WangY-2015},
also an example of hidden chaotic attractor in electromechanical model with no equilibria
was reported in a power system in 2001 \cite{Venkatasubramanian-2001}.

In this work we will consider three different systems which experience hidden oscillations in sense of mathematical definition of this term. But we will also show that some of these oscillations can be easily localized if physical nature of the process in such systems is taken into account.

\section{Translational oscillator--rotational actuator}

\begin{figure}[!h]
 \begin{center}                                                
   \includegraphics[width=0.4\textwidth]{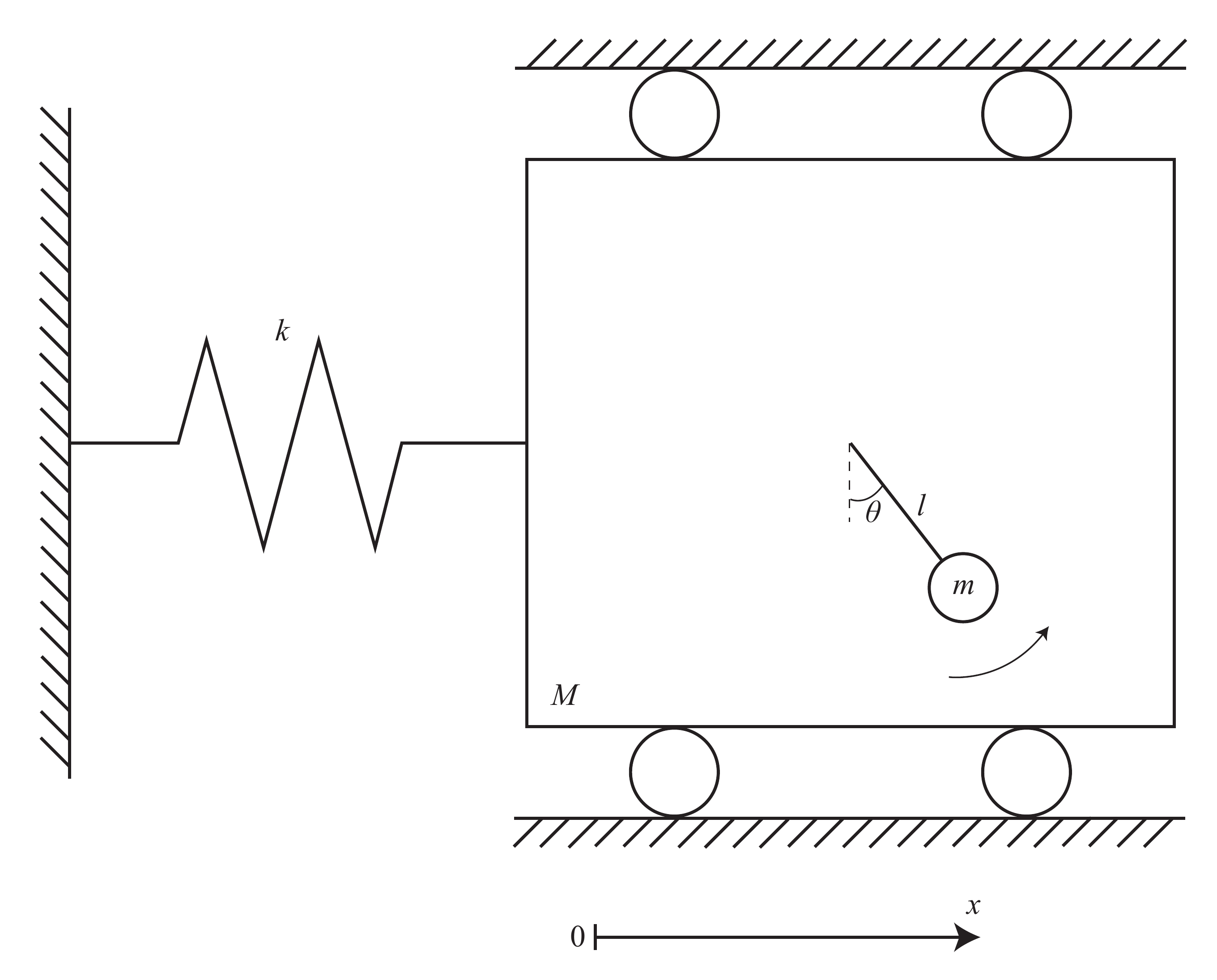}
 \end{center}
 \caption{Translational oscillator-rotational actuator scheme}
 \label{TORAscheme}
\end{figure}

Following the works \citep{EvanIwanowski-1976,FradkovTT-2011} let us consider the electromechanical system ``translational oscillator--rotational actuator'' (TORA) - see Fig.~\ref{TORAscheme}. It consists of DC motor which actuates the eccentric mass $m$ with eccentricity $l$ connected to the cart $M$. The cart is elastically connected to the wall with help of a string and moves only horizontally.
The equations of the system are as follows:
\begin{equation}\label{sys1}
\begin{array}{l}
(M+m) \ddot x + k_1 \dot x + ml(\ddot \theta \cos \theta - \dot \theta^2 \sin \theta) + kx = 0, \\
J \ddot \theta + k_\theta \dot \theta +ml \ddot x \cos \theta = u,
\end{array}
\end{equation}
Here $\theta$ is rotational angle of the rotor, $x$ is the displacement of the cart from its equilibrium position, $u$ is motor torque and $k$ is stiffness of the string, $k_1$ and $k_\theta$ are damping coefficients, $I$ is moment of inertia.

\begin{figure}[!h]
 \begin{center}                                                
   \includegraphics[width=0.4\textwidth]{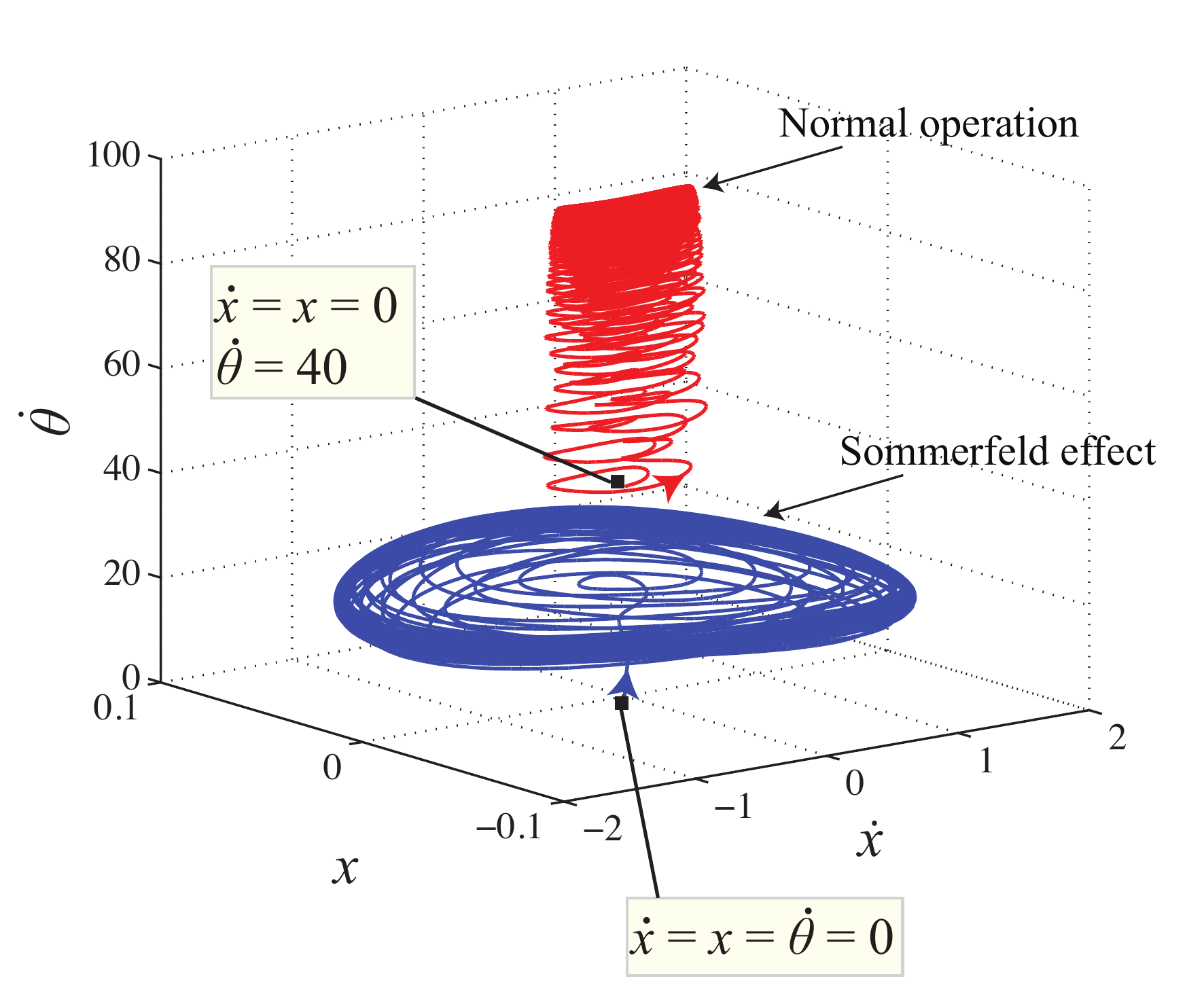}
 \end{center}
 \caption{Sommerfeld effect and normal operation in TORA}
 \label{SommNoSomm}
\end{figure}

Note that for $u \neq 0$ this system has no equilibria. Consider the following parameters of the system \citep{FradkovTT-2011}: $J = 0.014$, $M = 10.5$, $m_0 = 1.5$, $l = 0.04$, $k_\theta = 0.005$, $k =5300$, $k_1=5$. For $u = 0.48$ the Sommerfeld effect may be observed for initial data $\dot x =x = \theta =\dot \theta = 0$ (zero initial data correspond to typical start of the system, so it was easy to find this effect). But for other initial data $\dot x =x = \theta =  0, \dot \theta = 40$ we observe normal operation -- the achievement of desired rotational velocity of our mechanical system\footnote{The existence of both effects were reported in \citep{FradkovTT-2011}, but in our work the parameters are chosen more precisely}.
 In Fig.~\ref{SommNoSomm} the transient process for both initial data is shown, 
 in Fig.~\ref{SommNoSomm} we observe the attractors which are obtained after the transient process. Note that if we compare this situation to experiment of Sommerfeld we see that an effective level of about 1/4 (comparing to normal operation) is achieved here when Sommerfeld effect takes place.

\begin{figure}[!h]
 \begin{center}                                                
   \includegraphics[width=0.4\textwidth]{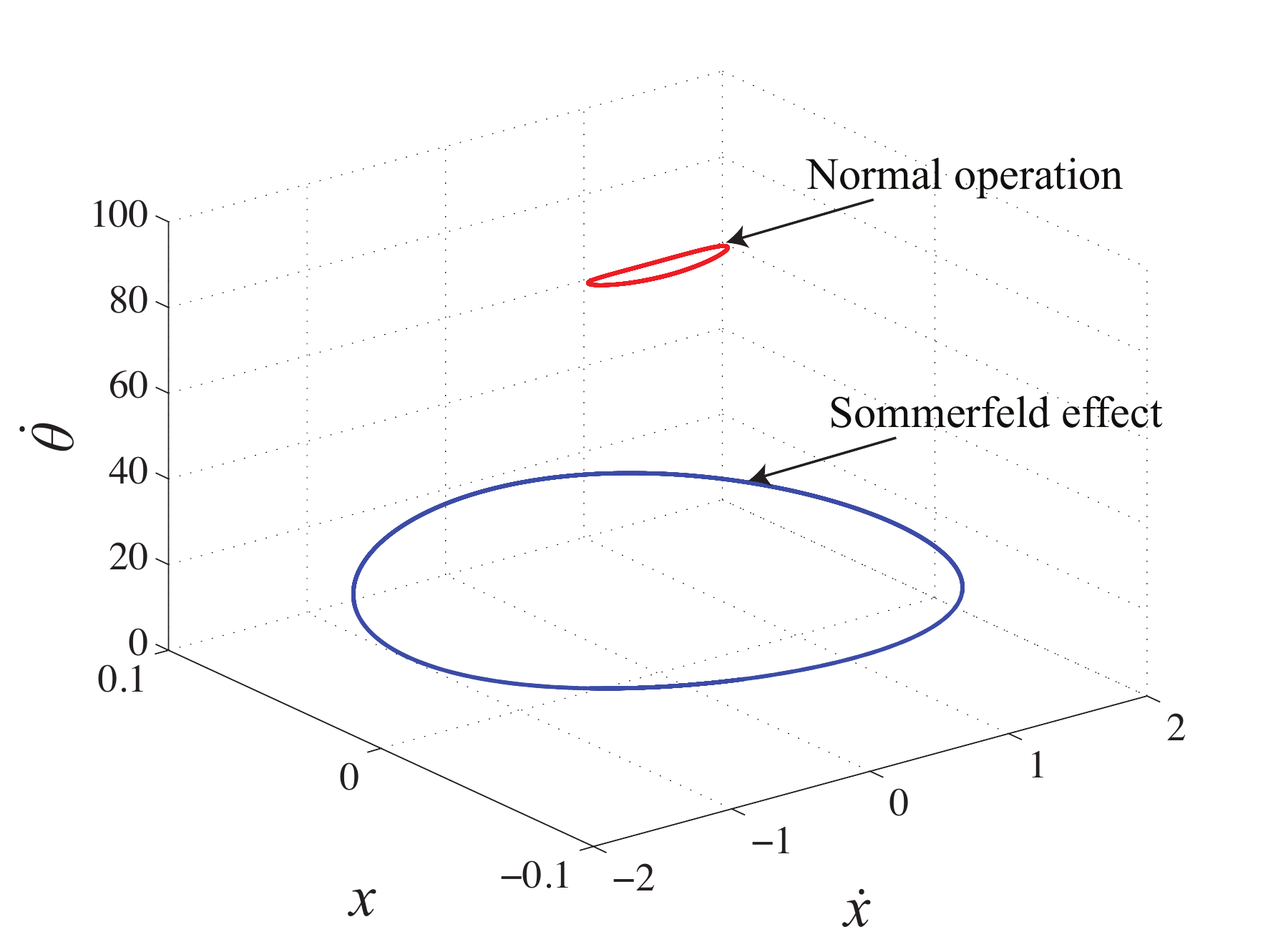}
 \end{center}
 \caption{Sommerfeld effect and normal operation in TORA -- localization after transient process}
 \label{SommNoSommLocal}
\end{figure}


\section{Drilling systems}
Let us consider now another electromechanical system -- drilling system. Drilling systems are widely used in oil and gas industry for drilling wells. The failures of drilling systems cause considerable time and expenditure loss for drilling companies, so the understanding of these failures is a very important task. Here we will consider two mathematical models of drilling systems and study their behaviour after operation start. For drilling systems two different ways of operation start are possible: no-load start and start with load. No-load start means that in initial moment of time there is no friction torque acting on the lower disc. Start with load implies start of the drilling with friction torque acting on the lower disc in initial moment of time (this case also corresponds to a sudden change of rock type). In contrast to TORA system these systems have stable equilibrium states which correspond to their normal operation.

\subsection{Drilling system actuated by DC motor}
\begin{figure}[!h]
 \begin{center}                                                
   \includegraphics[width=0.4\textwidth]{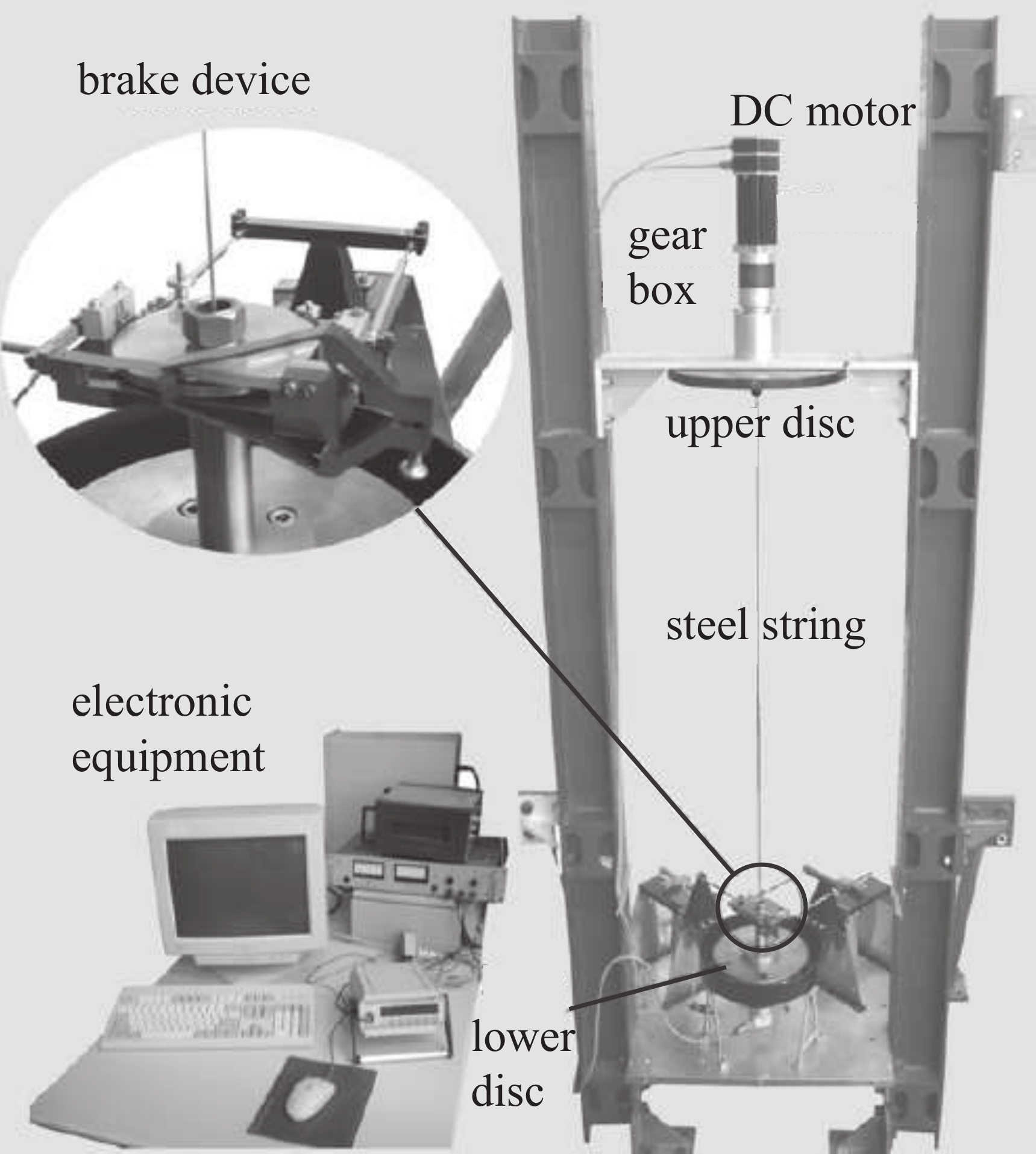}
 \end{center}
 \caption{Drilling system setup \citep{deBruin-2009}}
 \label{DrillingSetup}
\end{figure}

The first model of a drilling system was studied in \citep{deBruin-2009,Mihailovic-2004}. 
In these works the scientific group from Eindhoven University of Technology constructed and studied an experimental setup shown in Fig.~\ref{DrillingSetup}. The configuration of it can be recognised in the structure of drilling systems. The setup consists of two discs connected with a steel string. The upper disc represents rotary table of the drilling system and is actuated by a DC-motor. The lower disc represents bottom hole assembly. It is also assumed that the drill string in massless ad experiences only torsional deformation. The system is described by the following equations:
\begin{equation}\label{DrillingEind}
\begin{array}{l}
    		J_u \ddot \theta_u
 		       	+ k_\theta \left( \theta_u - \theta_l \right)
		        + b \left( \dot \theta_u - \dot \theta_l \right)
		        + T_{fu} \left( \dot \theta_u \right) - k_m v
		        = 0,
		\\
    		J_l\ddot \theta_l
       			- k_\theta \left( \theta_u - \theta_l \right)
        		- b \left( \dot \theta_u - \dot \theta_l \right)
        		+ T_{fl} \left( \dot \theta_l \right) = 0.
        \end{array}
\end{equation}
where $\theta_u(t)$ and $\theta_l(t)$ are angular displacements of the upper and lower discs with respect to the earth, $J_u$ and $J_l$ are constant inertia torques, $b$ is rotational friction, $k_\theta$  is the torsional spring stiffness, $k_m$ is the motor constant, $v$ is the constant input voltage. $T_{fu}$ and $T_{fl}$ are friction torques acting on the upper and on the lower disc respectively. Both friction torques $T_{fu}$ and $T_{fl}$ are obtained experimentally.
\begin{equation}
T_{fu}(\dot\theta_u)
			\in \left\{
				\begin{array}{ll}
					T_{cu}(\dot\theta_u){\rm sign}(\dot\theta_u),
						&  \dot\theta_u\neq 0 \\
					\left[-T_{su} +\Delta T_{su},T_{su} + \Delta T_{su}\right],
						& \dot\theta_u = 0,
				\end{array}
			\right.
\end{equation}
where
\begin{equation}
	T_{cu}(\dot\theta_u) =
				T_{su}+\Delta T_{su} {\rm sign} (\dot \theta_u) + b_u |\dot \theta_u| + \Delta b_u \dot \theta_u.	
\end{equation}
and
\begin{equation}
T_{fl}(\dot\theta_l)
			\in \left\{
				\begin{array}{ll}
					T_{cl}(\dot\theta_l){\rm sign}(\dot\theta_l),
						&  \dot\theta_l\neq 0 \\
					\left[-T_0,T_0\right],
						& \dot\theta_l = 0,
				\end{array}
			\right.
\end{equation}
where
\begin{equation}
	T_{cl}( \dot\theta_l) =
				\frac{T_0}{T_{sl}}(T_{pl}+(T_{sl}-T_{pl})
				e^{-|\frac{\dot\theta_l}{\omega_{sl}}|^{\delta_{sl}}}
				+ b_|\dot\theta_l).	
\end{equation}
Here $T_{su}$, $\Delta T_{su}$, $b_u$, $\Delta b_u$, $T_0$, $T{sl}$, $T_{pl}$, $\omega_{sl}$, $\delta_{sl}$, $b_l$ are constant parameters. Note that $T_{fu}$ and $T_{fl}$ are multi-valued functions, so special methods for numerical modelling of (\ref{DrillingEind}) are needed (see \citep{piiroinen2008event,Kiseleva2013J}).

\begin{figure}[!h]
 \begin{center}                                                
   \includegraphics[width=0.4\textwidth]{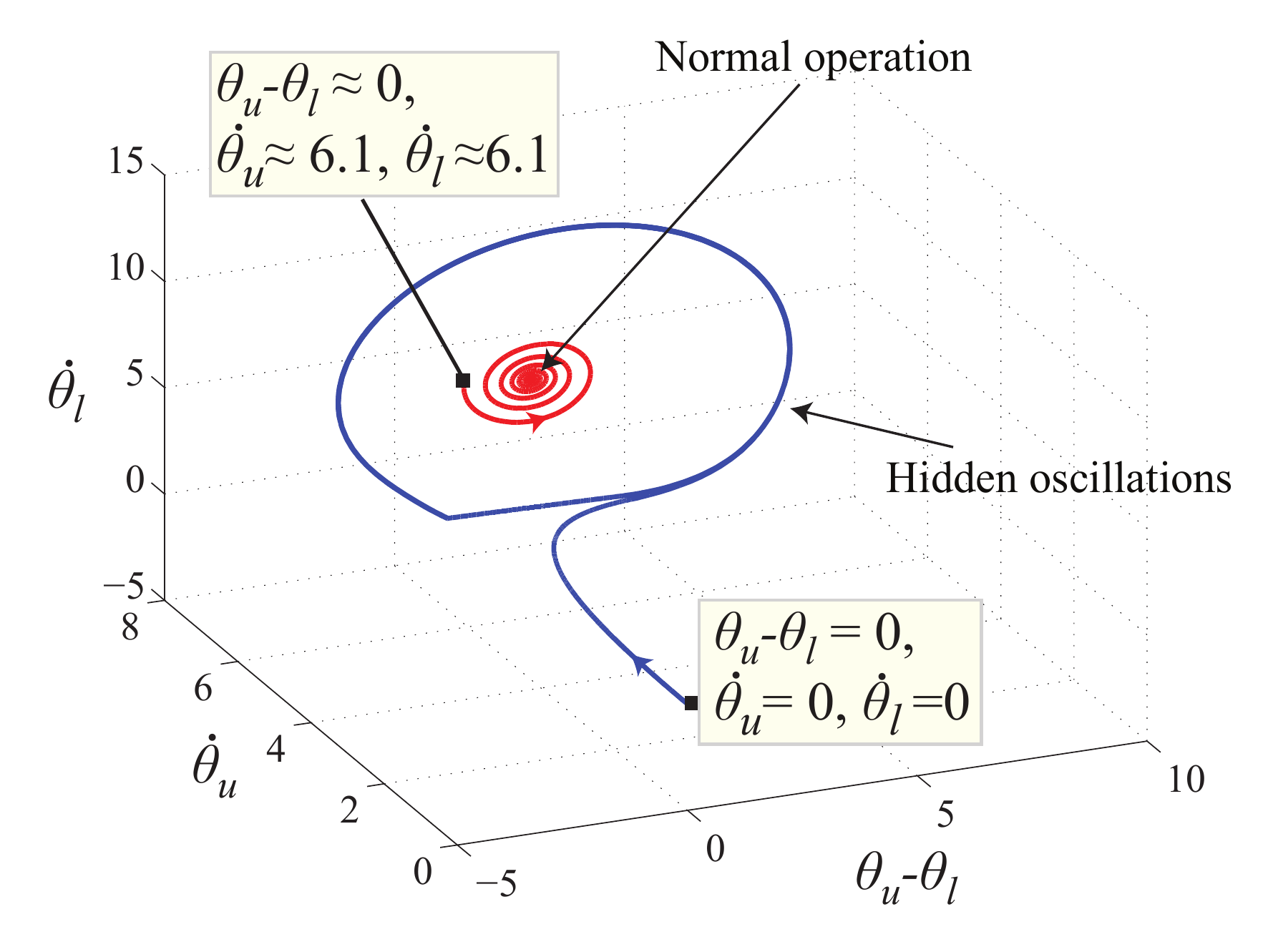}
 \end{center}
 \caption{Hidden oscillations and normal operation (corresponds to stable equilibrium state) in drilling system with DC motor}
 \label{DrillingEindhoven}
\end{figure}

Normal operation of the drilling system corresponds to rotation of both upper and lower discs with the same angular velocity with constant angular speed (i.e. the system reaches stable equilibrium state). Instead normal operation system may experience unwanted oscillations which lead to its failures. For modelling (\ref{DrillingEind}) we use the following parameters \citep{deBruin-2009}: $k_m = 4.3228$, $J_u = 0.4765$, $T_{su} = 0.37975$, $\Delta T_{su} = -0.00575$, $b_u = 2.4245$, $\Delta b_u = -0.0084$, $k_\theta = 0.075$, $b = 0$, $J_l = 0.035$, $T_{sl} = 0.26$, $T_{pl} = 0.05$, $\omega_{sl} = 2.2$, $\delta_{sl} = 1.5$, $b_l = 0.09$. For initial data $\theta_u-\theta_l = 0$, $\dot \theta_u = \dot \theta_l = 6.1$ (both upper and lower discs rotate with the same angular speed without angular displacement) after transient process the system enters normal operation mode (see Fig.~\ref{DrillingEindhoven}). But for same parameters and initial data $\theta_u-\theta_l = \dot \theta_u = \dot \theta_l = 0$ (initially discs don't rotate and there is no angular displacement between them) after transient process the system starts to experience stable hidden oscillations.

\subsection{Drilling system actuated by induction motor}
Consider now the modification of the drilling system studied above. Suppose it is driven by an induction motor (see schematic view of the system on Fig.~\ref{DrillingSystem}). In order to take into account the dynamics of the motor we write the equations in the following form \citep{KiselevaKKLS-2014,LeonovKKSZ-2014,Kiseleva2013J}:
\begin{equation}\label{DrillingInd}
\begin{array}{l}
			J_u \ddot\theta_u
		    	+ k_\theta \left(\theta_u - \theta_l \right)
  		      	+ b \left( \dot\theta_u - \dot\theta_l \right)
  		      \\	- nBS \sum\limits_{k=1}^{3} i_k
  		      		\sin \left(
		                \theta_u + \frac{2(k-1)\pi}{3}
		            \right) = 0,
		\cr
 		   	J_l \ddot{\theta_l}
 		   		- k_\theta \left( \theta_u - \theta_l \right)
  		      	- b \left( \dot\theta_u - \dot\theta_l \right)
 		       	+ T_{fl} \left( \omega + \dot{\theta_l} \right) = 0,
		\cr
		 	L \dot i_1 + (R+r) i_1 =
		 		- n B S \dot\theta_u \sin \theta_u,
		\cr
    		L \dot i_2 + (R+r) i_2 = - n B S \dot\theta_u
    			\sin \left(
           			 \theta_u + \frac{2\pi}{3}
        		\right),
		\cr
   			L \dot i_3 + (R+r) i_3 = - n B S \dot\theta_u
   				\sin \left(
           		 	\theta_u + \frac{4\pi}{3}
        		\right),
        \end{array}
\end{equation}
where $\theta_u(t)$ and $\theta_l(t)$ are angular displacements of the upper and lower discs with respect to the magnetic field rotating with constant speed $\omega = 2 \pi f/p$, where $f$ is the motor supply frequency, $p$ is the number of pairs of poles (usually not less than 8 pairs) \citep{2001-Springer-Leonhard}; $n$ is the number of turns in each coil; $B$ is an induction of magnetic field; $S$  is an area of one turn of coil; $i_k$ are currents in coils; $R$ is resistance of each coil; $r$ -- variable external resistance; $L$ -- inductance of each coil;
$J$ -- the moment of inertia of the rotor. Note that here angular displacements of the upper and lower discs with respect to the earth are $\theta_u(t)+\omega t$ and $\theta_l(t)+ \omega t$. Friction torque $T_{fl}$ acting on the lower disc is defined in the same way as in the previous model:
\begin{equation}
T_{fl}(\dot\theta_l+\omega)
			\in \left\{
				\begin{array}{ll}
					T_{cl}(\dot\theta_l+\omega){\rm sign}(\dot\theta_l+\omega),
						&  \dot\theta_l+\omega\neq 0 \\
					\left[-T_0,T_0\right],
						& \dot\theta_l+\omega = 0,
				\end{array}
			\right.
\end{equation}
where
\begin{equation}
	T_{cl}( \dot\theta_l+\omega) =
				\frac{T_0}{T_{sl}}(T_{pl}+(T_{sl}-T_{pl})
				e^{-|\frac{\dot\theta_l+\omega}{\omega_{sl}}|^{\delta_{sl}}}
				+ b_l(\dot\theta_l+\omega)).	
\end{equation}

\begin{figure}[!h]
 \begin{center}                                                
   \includegraphics[width=0.15\textwidth]{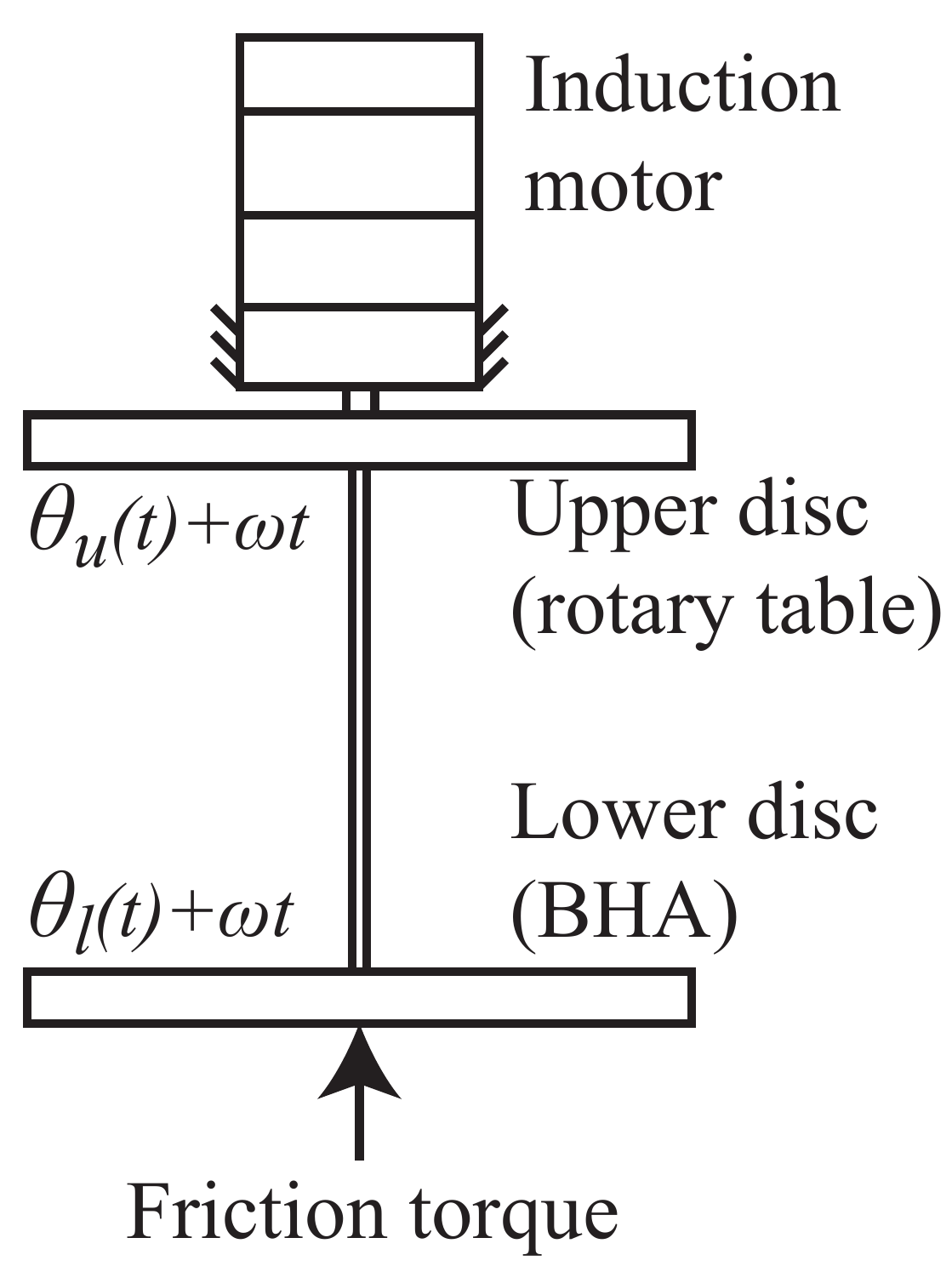}
 \end{center}
 \caption{Mathematical model of drilling system actuated by induction motor}
 \label{DrillingSystem}
\end{figure}

Let us model the system (\ref{DrillingInd}) with the following parameters: $T_0 = 0.25$, $c = 10$, $\omega = 8$, $J_u = 0.4765$, $J_l = 0.035$, $k = 0.075$, $a = 2.1$, $b = 0$, $T_{sl} = 0.26$, $T_{pl} = 0.05$, $\omega_{sl} = 2.2$, $\delta_{sl} = 1.5$, $b_l = 0.009$. For initial data $\theta$=$\theta_u-\theta_l$=0, $\omega_u = -\dot \theta_u = 0$ and $\omega_l = -\dot \theta_l = 0$ (rotation of both discs with the same speed without angular displacement) after the transient process the drilling system enters normal operation mode (see Fig.~\ref{DrillingSystem}). But for initial data $\theta$=0, $\omega_u = 8$ and $\omega_l = 8$ (initially discs don't rotate and there is no angular displacement) after the transient process the system starts to experience hidden oscillations, which may lead to break-down.

\begin{figure}[!h]
 \begin{center}                                                
   \includegraphics[width=0.4\textwidth]{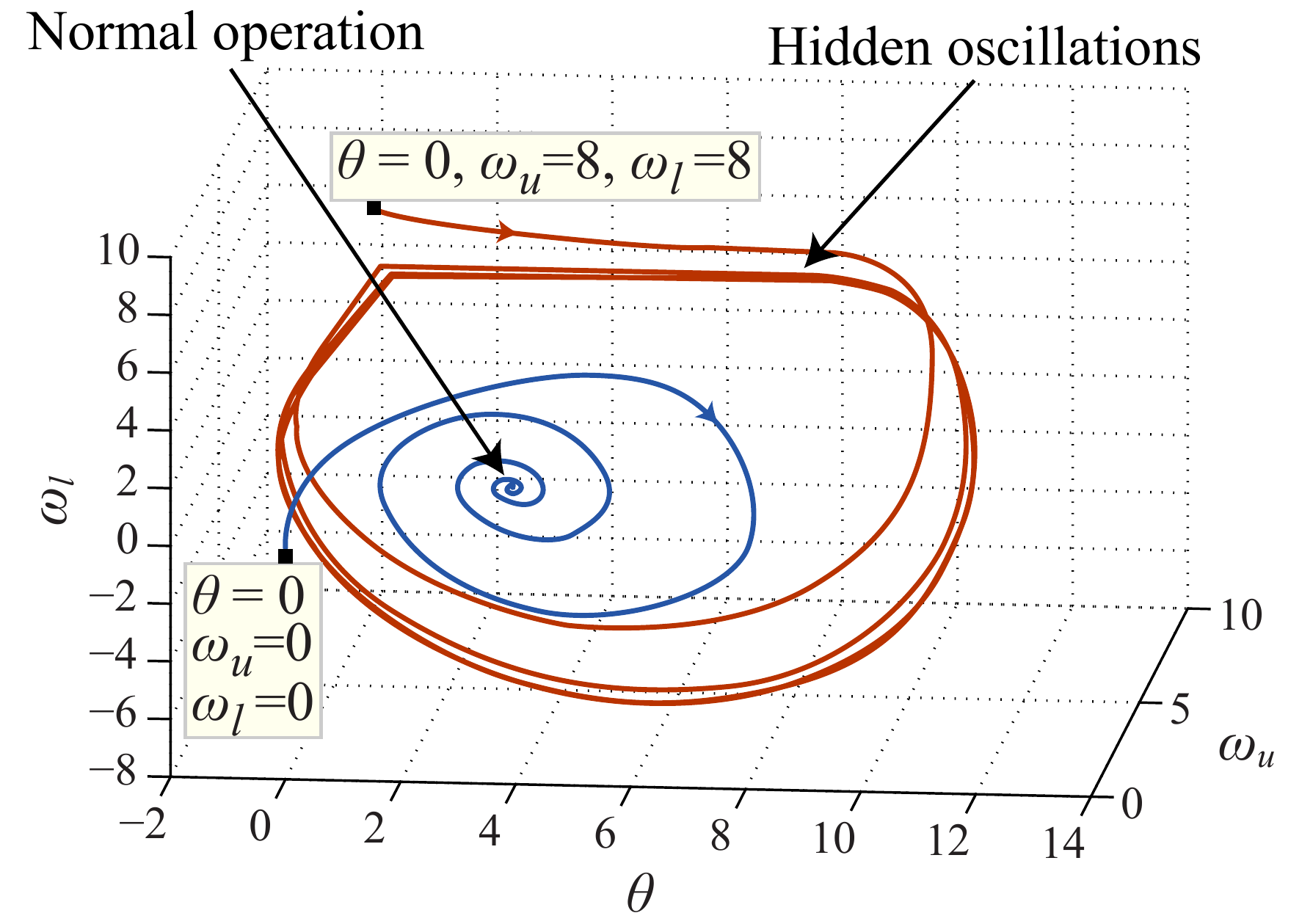}
 \end{center}
 \caption{Hidden oscillations and normal operation (corresponds to stable equilibrium state) in drilling system actuated by induction motor}
 \label{DrillingSystem}
\end{figure}

%

\section{Conclusions}
We modelled three different electromechanical systems. All of them have hidden oscillations in sense of mathematical definition. But from physical point of view some of these oscillations are easily localized, thus are not actually hidden. For example, for TORA system zero initial data correspond to typical start of the system, so Sommerfeld effect can be easily localized.
In our examples for drilling systems no-load start leads to normal operation and start with load (or change of rock type) leads to unwanted hidden oscillations.
Hence better understanding of physical nature of operation of electromechanical models may make it easier to find hidden attractors.

\bigskip
\section*{Acknowledgments}
 Authors were supported by Saint-Petersburg State University 
 and Russian Scientific Foundation. 


\end{document}